# Numerical Solutions of a Class of Second Order Boundary Value Problems on Using Bernoulli Polynomials

**Md. Shafiqul Islam[1*], Afroza Shirin [2]**

[1]Department of Mathematics, Dhaka University, Dhaka, Bangladesh.

[2]Department of Mathematics, Bangladesh University of Engineering and Technology, Dhaka, Bangladesh.

*Corresponding Author, Email: mdshafiqul@yahoo.com

## Abstract

The aim of this paper is to find the numerical solutions of the second order linear and nonlinear differential equations with *Dirichlet*, *Neumann* and *Robin* boundary conditions. We use the Bernoulli polynomials as linear combination to the approximate solutions of 2nd order boundary value problems. Here the Bernoulli polynomials over the interval [0, 1] are chosen as trial functions so that care has been taken to satisfy the corresponding homogeneous form of the *Dirichlet* boundary conditions in the Galerkin weighted residual method. In addition to that the given differential equation over arbitrary finite domain $[a, b]$ and the boundary conditions are converted into its equivalent form over the interval [0, 1]. All the formulas are verified by considering numerical examples. The approximate solutions are compared with the exact solutions, and also with the solutions of the existing methods. A reliable good accuracy is obtained in all cases.

## 1. Introduction

There are many linear and nonlinear problems in science and engineering, namely second order differential equations with various types of boundary conditions, are solved either analytically or numerically. In the literature of numerical analysis solving a two point second order boundary value problem (BVP) of differential equations, many authors have attempted to obtain higher accuracy rapidly by using a numerous methods. Among various numerical techniques, finite difference method has been widely used but it takes more computational costs to get high accuracy. In this method, a large number of parameters are required and it can not be used to evaluate the value of the desired points between two grid points. For this, Galerkin weighted residual method is widely used to find the approximate results to any point in the domain of the problem.

Since piecewise polynomials can be differentiated and integrated easily, and can be approximated any function to any accuracy desired [1], spline functions have been studied extensively in [2 – 9]. Solving BVP only with *Dirichlet* boundary conditions has been attempted in [4] while Bernstein polynomials [10, 11] have been used to solve the two point BVP very recently by the authors Bhatti and Bracken [1] rigorously by the Galerkin method. But it is limited only to second order BVP with *Dirichlet* boundary conditions, and to first order nonlinear differential equation. On the other hand, Ramadan *et al* [2] has studied linear BVP with *Neumann* boundary conditions using quadratic and cubic polynomial splines, and nonpolynomial splines. We have also found that the linear BVP with *Robin* (mixed) boundary conditions have been solved using finite difference method [12] and Sinc-Collocation method [13], respectively. Thus except [9], little attention has been given to solve the second order nonlinear BVP with *Dirichlet* and *Neumann* as well as *Robin* boundary conditions. Therefore, the purpose of this paper is to present the Galerkin weighted residual method to solve both linear and nonlinear second order BVP with all types of boundary conditions as well.

Besides spline functions and Bernstein polynomials, there is another type of piecewise continuous polynomials, namely Bernoulli polynomials, has been introduced by Atkinson in [14]. But none has attempted, to the knowledge of the present authors, using these polynomials to solve the second order BVP. Thus we concentrate in this paper rigorously to solve some linear and nonlinear BVP with various types of boundary conditions numerically, though it is originated in [1].

However, in this paper, we first give an introduction of Bernoulli polynomials, and then we formulate the Galerkin approximation method using Bernoulli polynomials. We derive the individual formulas for each BVP consisting of *Dirichlet*, *Neumann* and *Robin* boundary conditions, respectively. Numerical examples, for both linear and nonlinear boundary value problems, are considered to verify the effectiveness of the derived formulas, and are also compared with the exact solutions. All the computations are performed by *MATHEMATICA*.

## 2. Bernoulli Polynomials

The Bernoulli polynomials [14, p. 284] of degree $n$ can be defined over the interval [0, 1] implicitly by

$$Br_n(x) = \sum_{k=0}^{n} \binom{n}{k} b_k x^{n-k} \tag{1a}$$

where, $b_k$ are Bernoulli numbers given by

$$b_0 = 1 \quad \text{and} \quad b_k = -\int_0^1 Br_k(x)\,dx \qquad k \geq 1. \tag{1b}$$

Also Equation (1) can be written explicitly as

$$Br_0(x) = 1$$

$$Br_m(x) = \sum_{n=0}^{m} \frac{1}{n+1} \sum_{k=0}^{n} (-1)^k \binom{n}{k} (x+k)^m - \sum_{n=0}^{m} \frac{1}{n+1} \sum_{k=0}^{n} (-1)^k \binom{n}{k} k^m, m \geq 1 \tag{2}$$

The the first 11 Bernoulli polynomials are given bellow:

$$Br_0(x) = 1 \qquad Br_5(x) = -\frac{x}{6} + \frac{5x^3}{3} - \frac{5x^4}{2} + x^5$$

$$Br_1(x) = x \qquad Br_6(x) = -\frac{x^2}{2} + \frac{5x^4}{2} - 3x^5 + x^6$$

$$Br_2(x) = -x + x^2 \qquad Br_7(x) = \frac{x}{6} - \frac{7x^3}{6} + \frac{7x^5}{2} - \frac{7x^6}{2} + x^7$$

$$Br_3(x) = \frac{x}{2} - \frac{3x^2}{2} + x^3 \qquad Br_8(x) = \frac{2x^2}{3} - \frac{7x^4}{3} + \frac{14x^6}{3} - 4x^7 + x^8$$

$$Br_4(x) = x^2 - 2x^3 + x^4 \qquad Br_9(x) = -\frac{3x}{10} + 2x^3 - \frac{21x^5}{5} + 6x^7 - \frac{9x^8}{2} + x^9$$

$$Br_{10}(x) = -\frac{3x^2}{2} + 5x^4 - 7x^6 + \frac{15x^8}{2} - 5x^9 + x^{10}$$

Since Bernoulli polynomials have special properties at $x = 0$ and $x = 1$:
$Br_n(0) = 0,\ n \geq 1$, and $Br_n(1) = 0,\ n \geq 2$ respectively, so that they can be used as a set of basis functions to satisfy the corresponding homogeneous form of the *Dirichlet* boundary conditions to derive the matrix formulation of second order BVP over the interval [0, 1].

## 3. Formulation of second order BVP

We consider the general second order linear BVP [15]:

$$-\frac{d}{dx}\left(p(x)\frac{du}{dx}\right) + q(x)u = r(x), \qquad a < x < b, \tag{3a}$$

$$\alpha_0 u(a) + \alpha_1 u'(a) = c_1, \quad \beta_0 u(b) + \beta_1 u'(b) = c_2, \tag{3b}$$

where $p(x)$, $q(x)$ and $r(x)$ are specified continuous functions, and $\alpha_0$, $\alpha_1$, $\beta_0$, $\beta_1$, $c_1$, $c_2$ are specified numbers. Since our aim is to use the Bernoulli polynomials as trial functions which are derived over the interval [0, 1], so the BVP (3) is to be converted to an equivalent problem on [0, 1] by replacing $x$ by $(b-a)x + a$, and thus we have:

$$-\frac{d}{dx}\left(\tilde{p}(x)\frac{du}{dx}\right) + \tilde{q}(x)u = \tilde{r}(x), \qquad 0 < x < 1 \tag{4a}$$

$$\alpha_0 u(0) + \frac{\alpha_1}{b-a}u'(0) = c_1, \qquad \beta_0 u(1) + \frac{\beta_1}{b-a}u'(1) = c_2 \tag{4b}$$

where,

$$\tilde{p}(x) = \frac{1}{(b-a)^2}\, p\big((b-a)x + a\big), \quad \tilde{q}(x) = q\big((b-a)x + a\big), \quad \tilde{r}(x) = r\big((b-a)x + a\big) \tag{4c}$$

Using Bernoulli polynomials, $Br_i(x)$ in Equation (2), we assume an approximate solution in a form,

$$\tilde{u}(x) = \sum_{i=0}^{n} a_i Br_i(x), \qquad n \geq 1, \tag{5}$$

Now the weighted residual equations corresponding to the differential equation (4a) given by

$$\int_0^1 \left[-\frac{d}{dx}\left(\tilde{p}(x)\frac{d\tilde{u}}{dx}\right) + \tilde{q}(x)\tilde{u} - \tilde{r}(x)\right] Br_j(x)\,dx = 0, \quad j = 1,2,3,\cdots,n \tag{6}$$

Since from (5), we have

$$\frac{d\tilde{u}}{dx} = \sum_{i=0}^{n} a_i \frac{dBr_i}{dx}, \quad \tilde{u}(0) = \sum_{i=0}^{n} a_i Br_i(0) \text{ and } \tilde{u}(1) = \sum_{i=0}^{n} a_i Br_i(1)$$

After minor simplification, from (6) we can obtain

$$\sum_{i=0}^{n}\left[\int_0^1\left[\widetilde{p}(x)\frac{dBr_i}{dx}\frac{dBr_j}{dx}+\widetilde{q}(x)Br_i(x)Br_j(x)\right]dx + \frac{\beta_0(b-a)\widetilde{p}(1)Br_i(1)Br_j(1)}{\beta_1}-\frac{\alpha_0(b-a)\widetilde{p}(0)Br_i(0)Br_j(0)}{\alpha_1}\right]a_i$$

$$=\int_0^1\widetilde{r}(x)Br_j(x)dx+\frac{c_2(b-a)\widetilde{p}(1)Br_j(1)}{\beta_1}-\frac{c_1(b-a)\widetilde{p}(0)Br_j(0)}{\alpha_1} \qquad (7)$$

or, equivalently in matrix form,

$$\sum_{i=0}^{n}D_{i,j}\ a_i=F_j,\ \ j=0,1,2,...,n \qquad (8a)$$

where,

$$D_{i,j}=\int_0^1\left[\widetilde{p}(x)\frac{dBr_i}{dx}\frac{dBr_j}{dx}+\widetilde{q}(x)Br_i(x)Br_j(x)\right]dx+\frac{\beta_0(b-a)\widetilde{p}(1)Br_i(1)Br_j(1)}{\beta_1}$$

$$-\frac{\alpha_0(b-a)\widetilde{p}(0)Br_i(0)Br_j(0)}{\alpha_1}\quad i,j=0,1,2,....,n \qquad (8b)$$

$$F_j=\int_0^1\widetilde{r}(x)Br_j(x)dx+\frac{c_2(b-a)\widetilde{p}(1)Br_j(1)}{\beta_1}-\frac{c_1(b-a)\widetilde{p}(0)Br_j(0)}{\alpha_1},\quad j=0,1,2,....,n \qquad (8c)$$

Solving the system (8a), we find the values of the parameters $a_i$ $(i=0,1,2,....,n)$, and then substituting these parameters in (5), we get the approximate solution of the BVP (4). If we replace $x$ by $\frac{x-a}{b-a}$ in $\widetilde{u}(x)$, then we get the desired approximate solution of the BVP (3).

The absolute error, $E$ of this formulation is defined by

$$E=\left|u(x)-\widetilde{u}(x)\right|.$$

Now we discuss the various types of BVP using different boundary conditions as follows:

**Case 1:** The matrix formulation with the *Robin* (mixed) boundary conditions**,** (*i.e.*, $\alpha_0\neq 0$, $\alpha_1\neq 0$, $\beta_0\neq 0$, $\beta_1\neq 0$), are already defined in Equations (8).

**Case 2:** The matrix formulation of the differential equation (3a) with the *Dirichlet* boundary conditions (*i.e.*, $\alpha_0\neq 0$, $\alpha_1\neq 0$, $\beta_0\neq 0$, $\beta_1\neq 0$), is given by

$$\sum_{i=2}^{n}D_{i,j}\ a_i=F_j,\qquad j=2,3,.....,n \qquad (9a)$$

where,

$$D_{i,j}=\int_0^1\left[\widetilde{p}(x)\frac{dBr_i}{dx}\frac{dBr_j}{dx}+\widetilde{q}(x)Br_i(x)Br_j(x)\right]dx,\qquad i,j=2,3,....,n. \qquad (9b)$$

$$F_j=\int_0^1\widetilde{r}(x)Br_j(x)dx-\int_0^1\left[\widetilde{p}(x)\frac{d\theta_0}{dx}\frac{dBr_j}{dx}+\widetilde{q}(x)\theta_0(x)Br_j(x)\right]dx,\ \ j=2,3,....,n \qquad (9c)$$

**Case 3:** The approximate solution of the differential equation (3a) consisting of *Neumann* boundary conditions (*i.e.,* $\alpha_0 = 0$, $\alpha_1 \neq 0$, $\beta_0 = 0$, $\beta_1 \neq 0$) can be obtained by putting $\alpha_0 = 0$ and $\beta_0 = 0$ in the Equation (8) as

$$\sum_{i=0}^{n} D_{i,j} \; a_i = F_j, \; j = 0,1,2,....,n \tag{10a}$$

where

$$D_{i,j} = \int_0^1 \left[ \tilde{p}(x) \frac{dBr_i}{dx} \frac{dBr_j}{dx} + \tilde{q}(x) Br_i(x) Br_j(x) \right] dx, \quad i,j = 0,1,2,.....,n \tag{10b}$$

$$F_j = \int_0^1 \tilde{r}(x) Br_j(x) dx + \frac{c_2 (b-a) \tilde{p}(1) Br_j(1)}{\beta_1} - \frac{c_1 (b-a) \tilde{p}(0) Br_j(0)}{\alpha_1}, \quad j = 0,1,2,....,n \tag{10c}$$

Similar formulation for nonlinear BVP using the Bernoulli polynomials can be derived, which will be discussed through numerical examples in the next section.

## 4. Numerical Examples

In this section, we explain four linear and two nonlinear boundary value problems which are available in the existing literatures, considering three types of boundary conditions to verify the effectiveness of the present formulations described in the previous sections. The convergence of each linear BVP is calculated by

$$E = \left| u_{n+1}(x) - u_n(x) \right| < \delta,$$

where $u_n(x)$ denotes the approximate solution by the proposed method using $n$-th degree polynomial approximation. The convergence of nonlinear BVP is assumed when the absolute error of two consecutive iterations is recorded below the convergence criterion $\delta$ such that

$$\left| \tilde{u}^{N+1} - \tilde{u}^{N} \right| < \delta$$

where $N$ is the Newton's iteration number and $\delta$ varies from $10^{-8}$.

**Example 1.** First we consider the BVP with *Robin* boundary conditions [15]:

$$-\frac{d^2 u}{dx^2} + u = 2\cos x, \qquad \pi/2 < x < \pi \tag{11a}$$

$$u'(\pi/2) + 3u(\pi/2) = -1, \qquad u'(\pi) + 4u(\pi) = -4, \tag{11b}$$

and the exact solution is $u(x) = \cos x$.

The BVP (11) over [0, 1] is equivalently to the BVP

$$-\frac{1}{(\pi/2)^2} \frac{d^2 u}{dx^2} + u = 2\cos(\frac{\pi}{2} x + \frac{\pi}{2}), 0 < x < 1 \tag{12a}$$

$$\frac{2}{\pi} u'(0) + 3u(0) = -1, \qquad \frac{2}{\pi} u'(1) + 4u(1) = -4 \tag{12b}$$

Using the formula derived in equations (8) and using different number of Bernoulli polynomials, the approximate solutions are summarized in **Table 1**. It is observed that the accuracy is found nearly the order $10^{-6}, 10^{-9}$ and $10^{-11}$ on using 5, 7 and 9 Bernoulli polynomials, respectively.

**Table 1. Approximate solutions and absolute differences for the example 1.**

| $x$ | Approximate | Absolute Error | Approximate | Absolute Error | Approximate | Absolute Error |
|---|---|---|---|---|---|---|
| | **Bernoulli polynomials, 5** | | **Bernoulli polynomials, 7** | | **Bernoulli polynomials, 9** | |
| π/2 | 00.0000000000 | 0.0000000000 | 00.000000000 | 0.0000000000 | 0.0000000000 | 0.00000000000 |
| 11π/20 | -0.1564317160 | $2.749065\times10^{-6}$ | -0.1564344475 | $1.755782\times10^{-8}$ | -0.1564344650 | $4.925213\times10^{-11}$ |
| 3π/5 | -0.3090255130 | $8.518602\times10^{-6}$ | -0.3090169914 | $3.001562\times10^{-9}$ | -0.3090169944 | $5.283479\times10^{-11}$ |
| 13π/20 | -0.4539971315 | $6.631794\times10^{-6}$ | -0.4539905271 | $2.730026\times10^{-8}$ | -0.4539904998 | $3.642830\times10^{-11}$ |
| 7π/10 | -0.5877810347 | $4.217626\times10^{-6}$ | -0.587785252 | $7.19809\times10^{-10}$ | -0.5877852522 | $1.148404\times10^{-10}$ |
| 3π/4 | -0.7070961552 | 0.0000110000 | -0.7071067522 | $2.896640\times10^{-8}$ | -0.7071067812 | $9.037660\times10^{-12}$ |
| 4π/5 | -0.8090116456 | $5.348787\times10^{-6}$ | -0.8090169912 | $3.186240\times10^{-9}$ | -0.8090169945 | $1.525206\times10^{-10}$ |
| 17π/20 | -0.8910126268 | $6.102617\times10^{-6}$ | -0.8910065523 | $2.810947\times10^{-8}$ | -0.8910065241 | $4.453860\times10^{-11}$ |
| 9π/10 | -0.9510659385 | $9.422186\times10^{-6}$ | -0.9510565150 | $1.320862\times10^{-9}$ | -0.9510565162 | $1.337538\times10^{-10}$ |
| 19π/20 | -0.9876858881 | $2.452495\times10^{-6}$ | -0.9876883211 | $1.952363\times10^{-8}$ | -0.9876883408 | $1.614215\times10^{-10}$ |
| π | -1.0000000000 | 0.0000000000 | -1.0000000000 | 0.0000000000 | -1.0000000000 | 0.00000000000 |

**Example 2.** We consider the BVP with *Dirichlet* boundary conditions [1]:

$$\frac{d^2u}{dx^2} + u = x^2 e^{-x}, \qquad 0 < x < 10 \tag{13a}$$

$$u(0) = 0, \qquad u(10) = 0 \tag{13b}$$

The exact solution is:

$$\frac{1}{2}e^{-x}\left[1 + 2x + x^2 - e^x \cos x - 2e^x \left\{-\frac{\cot 10}{2} + \frac{121 \operatorname{cosec} 10}{2e^{10}}\right\} \sin x\right]$$

The BVP (13) is equivalent to the BVP

$$\frac{1}{100}\frac{d^2u}{dx^2} + u = 100x^2 e^{-10x}, \qquad 0 < x < 1 \tag{14a}$$

$$u(0) = u(1) = 0, \tag{14b}$$

Using the formula derived in Equations (9), the approximate solutions, shown in **Table 2**, are obtained using 8, 10 and 15 Bernoulli polynomials with accuracy up to 3, 4 and 6 significant digits, respectively. It is observed that using 21 Bernstein polynomials, the accuracy is found nearly the order of $10^{-5}$ in [1].

**Example 3.** In this case **w**e consider the BVP with *Dirichlet* boundary conditions [4]:

$$\frac{d^2u}{dx^2} = \frac{2}{x^2}u - \frac{1}{x}, \qquad 2 < x < 3 \tag{15a}$$

$$u(2) = 0, \qquad u(3) = 0 \tag{15b}$$

The exact solution is:

$$u(x)=\frac{1}{38}\left[-5x^2+19x-\frac{36}{x}\right]$$

**Table 2. Approximate solutions and absolute differences for the example 2.**

| $x$ | **Approximate** | **Absolute Error** | **Approximate** | **Absolute Error** | **Approximate** | **Absolute Error** |
|---|---|---|---|---|---|---|
| | **Bernoulli polynomials, 8** | | **Bernoulli polynomials, 10** | | **Bernoulli polynomials, 15** | |
| 0. | 0.0000000000 | 0.0000000000 | 0.0000000000 | 0.0000000000 | 0.0000000000 | 0.0000000000 |
| 1. | 1.1136094978 | 0.0051685418 | 1.1187686211 | $9.418441\times10^{-6}$ | 1.1187829850 | $4.945487\times10^{-6}$ |
| 2. | 1.5330008942 | 0.0100998520 | 1.5229276754 | $2.663320\times10^{-5}$ | 1.5228972291 | $3.813112\times10^{-6}$ |
| 3. | 1.0001516651 | 0.0026819237 | 1.0028553869 | $2.179813\times10^{-5}$ | 1.0028378028 | $4.213962\times10^{-6}$ |
| 4. | -0.0413011761 | 0.0096199100 | -0.0317862308 | $1.049646\times10^{-4}$ | -0.0316869594 | $5.693175\times10^{-6}$ |
| 5. | -0.7601267725 | 0.0047616938 | -0.7647484696 | $1.399968\times10^{-4}$ | -0.7648818892 | $6.577111\times10^{-6}$ |
| 6. | -0.6303887189 | 0.0058561235 | -0.6363262078 | $8.136546\times10^{-5}$ | -0.6362508546 | $6.012205\times10^{-6}$ |
| 7. | 0.1575008515 | 0.0046972725 | 0.1621779732 | $2.015082\times10^{-5}$ | 0.1622028015 | $4.677455\times10^{-6}$ |
| 8. | 0.8534716130 | 0.0008286515 | 0.8543763192 | $7.605468\times10^{-5}$ | 0.8542961563 | $4.108267\times10^{-6}$ |
| 9. | 0.7833934560 | 0.0017614110 | 0.7815841935 | $4.785153\times10^{-5}$ | 0.7816365818 | $4.536826\times10^{-6}$ |
| 10. | 0.0000000000 | 0.0000000000 | 0.0000000000 | 0.0000000000 | 0.0000000000 | 0.0000000000 |

The BVP (15) is equivalent to the BVP

$$\frac{d^2u}{dx^2}=\frac{2}{(x+2)^2}u-\frac{1}{x+2},\qquad 0<x<1, \tag{16a}$$

$$u(0)=0,\qquad u(1)=0 \tag{16b}$$

Using the formula derived in Equations (9), the approximate solutions, shown in **Table 3**, are obtained on using 5, 7 and 10 Bernoulli polynomials, and the accuracy is observed nearly 7, 8 and 9 decimal places, respectively. On the contrary, the error is obtained nearly $10^{-10}$ by Arshad [10] with $h=1/32$, where $h=(b-a)/N$, $a$ and $b$ are the endpoints of the domain and $N$ is the number of subdivision of intervals [*a, b*].

**Example 4.** We consider the BVP with *Neumann* boundary conditions [2]:

$$\frac{d^2u}{dx^2}+u=-1,\quad 0\le x\le 1 \tag{17a}$$

$$u'(0)=\frac{1-\cos 1}{\sin 1},\qquad u'(1)=-\frac{1-\cos 1}{\sin 1} \tag{17b}$$

whose exact solution is, $u(x)=\cos x+\dfrac{1-\cos 1}{\sin 1}\sin x-1$.

Using the formula given in Equations (10), the approximate solutions, shown in **Table 4**, are obtained on using 5, 7 and 10 Bernoulli polynomials with the remarkable accuracy nearly the order of $10^{-7}$, $10^{-10}$ and $10^{-14}$. On the other hand, Ramadan *et al* [6] has found nearly the accuracy of order $10^{-6}$ *and* $10^{-6}$, and $10^{-8}$ on using quadratic and cubic polynomial splines, and nonpolynomial spline, respectively with $h = 1/128$, where $h = (b-a)/N$, $a$ and $b$ are the endpoints of the domain and $N$ is the number of subdivision of intervals [$a$, $b$].

**Table 3. Approximate solutions and absolute differences for the example 3.**

| $x$ | Approximate | Absolute Error | Approximate | Absolute Error | Approximate | Absolute Error |
|---|---|---|---|---|---|---|
| | Bernoulli polynomials = 5; | | Bernoulli polynomials = 7; | | Bernoulli polynomials = 1 0; | |
| 2.0 | 0.0000000000 | 0.0000000000 | 0.0000000000 | 0.0000000000 | 0.0000000000 | 0.0000000000 |
| 2.1 | 0.0186087702 | $2.523743\times10^{-7}$ | 0.0186090317 | $9.107887\times10^{-9}$ | 0.0186090279 | $5.349871\times10^{-9}$ |
| 2.2 | 0.0325370415 | $1.156379\times10^{-6}$ | 0.0325358850 | $1.414425\times10^{-10}$ | 0.0325358805 | $4.662072\times10^{-9}$ |
| 2.3 | 0.0420487288 | $6.738722\times10^{-7}$ | 0.0420480426 | $1.229932\times10^{-8}$ | 0.0420480538 | $1.138989\times10^{-9}$ |
| 2.4 | 0.0473677309 | $6.901458\times10^{-7}$ | 0.0473684233 | $2.207156\times10^{-9}$ | 0.0473684265 | $5.421880\times10^{-9}$ |
| 2.5 | 0.0486829640 | $1.246501\times10^{-6}$ | 0.0486842230 | $1.246639\times10^{-8}$ | 0.0486842096 | $9.479982\times10^{-10}$ |
| 2.6 | 0.0461533943 | $4.518647\times10^{-7}$ | 0.0461538458 | $3.633628\times10^{-10}$ | 0.0461538418 | $4.318545\times10^{-9}$ |
| 2.7 | 0.0399130707 | $7.900046\times10^{-7}$ | 0.0399122691 | $1.157579\times10^{-8}$ | 0.0399122828 | $2.126246\times10^{-9}$ |
| 2.8 | 0.0300761580 | $9.700405\times10^{-7}$ | 0.0300751896 | $1.649355\times10^{-9}$ | 0.0300751900 | $2.000938\times10^{-9}$ |
| 2.9 | 0.0167419695 | $3.172309\times10^{-7}$ | 0.0167422939 | $7.169105\times10^{-9}$ | 0.0167422842 | $2.571030\times10^{-9}$ |
| 3.0 | 0.0000000000 | 0.0000000000 | 0.0000000000 | 0.0000000000 | 0.0000000000 | 0.0000000000 |

**Table 4. Approximate solutions and absolute differences for the example 4.**

| $x$ | Approximate | Absolute Error | Approximate | Absolute Error | Approximate | Absolute Error |
|---|---|---|---|---|---|---|
| | Bernoulli polynomials = 5 | | Bernoulli polynomial = 7 | | Bernoulli polynomials = 10 | |
| 0.0 | 0.0000000000 | $1.153937\times10^{-10}$ | 0.0000000000 | 0.0000000000 | 0.0000000000 | 0.0000000000 |
| 0.1 | 0.0495431439 | $2.654280\times10^{-7}$ | 0.0495434086 | $8.177960\times10^{-10}$ | 0.0495434094 | $1.932482\times10^{-14}$ |
| 0.2 | 0.0886011150 | $9.870426\times10^{-7}$ | 0.0886001278 | $8.626597\times10^{-11}$ | 0.0886001279 | $2.672862\times10^{-14}$ |
| 0.3 | 0.1167806021 | $6.883117\times10^{-7}$ | 0.1167799150 | $1.222063\times10^{-9}$ | 0.1167799138 | $1.992850\times10^{-14}$ |
| 0.4 | 0.1338006690 | $5.349698\times10^{-7}$ | 0.1338012039 | $9.031828\times10^{-11}$ | 0.1338012040 | $1.013079\times10^{-14}$ |
| 0.5 | 0.1394927538 | $1.173569\times10^{-6}$ | 0.1394939260 | $1.280755\times10^{-9}$ | 0.1394939273 | $2.942091\times10^{-14}$ |
| 0.6 | 0.1338006690 | $5.349698\times10^{-7}$ | 0.1338012039 | $9.031831\times10^{-11}$ | 0.1338012040 | $1.010303\times10^{-14}$ |
| 0.7 | 0.1167806021 | $6.883117\times10^{-7}$ | 0.1167799150 | $1.222063\times10^{-9}$ | 0.1167799138 | $1.981748\times10^{-14}$ |
| 0.8 | 0.0886011150 | $9.870426\times10^{-7}$ | 0.0886001278 | $8.626594\times10^{-11}$ | 0.0886001279 | $2.670086\times10^{-14}$ |
| 0.9 | 0.0495431439 | $2.654280\times10^{-7}$ | 0.0495434085 | $8.177957\times10^{-10}$ | 0.0495434094 | $1.942890\times10^{-14}$ |
| 1.0 | 0.0000000000 | 0.0000000000 | 0.0000000000 | 0.0000000000 | 0.0000000000 | 0.0000000000 |

We now also implement the procedure described in section 3 to find the numerical solutions of two nonlinear second order boundary value problems.

**Example 5.** We consider a nonlinear BVP with *Dirichlet* boundary conditions [16]

$$\frac{d^2u}{dx^2} + \frac{1}{8}u\frac{du}{dx} = (4 + \frac{1}{4}x^3), \quad 1 < x < 3 \tag{18a}$$

$$u(1) = 17 \qquad \text{and} \qquad u(3) = 43/3 \tag{18b}$$

The exact solution of the problem is given by

$$u(x) = x^2 + \frac{16}{x}$$

To use Bernoulli polynomials, first we convert the BVP (18) to an equivalent BVP on $[0, 1]$ by replacing $x$ by $2x+1$ such that

$$\frac{d^2u}{dx^2} + \frac{1}{4}u\frac{du}{dx} = 16 + (2x+1)^3, \qquad 0 < x < 1, \tag{19a}$$

$$u(0) = 17 \qquad \text{and} \qquad u(1) = 43/3\,. \tag{19b}$$

Assume that the approximate solution of (19) using Bernoulli polynomials is given by

$$\tilde{u}(x) = \theta_0(x) + \sum_{i=2}^{n} a_i Br_i(x), \qquad n \geq 2, \tag{20}$$

where $\theta_0(x) = 17 - 8x/3$ is specified by the *Dirichlet* boundary conditions at $x = 0$ and $x = 1$, and $Br_i(0) = Br_i(1) = 0$ for each $i = 2,3,\ldots,n.$

The weighted residual equations of (19*a*) corresponding to the approximation (20), given by

$$\int_0^1 \left( \frac{d^2\tilde{u}}{dx^2} + \frac{1}{4}\tilde{u}\frac{d\tilde{u}}{dx} - [16 + (2x+1)^3] \right) Br_k(x)dx = 0\,, \qquad k = 2,3,\ldots,n. \tag{21}$$

Exploiting integration by parts with minor simplifications, we obtain

$$\sum_{i=2}^{n}\left[\int_0^1\left[-\frac{dBr_i}{dx}\frac{dBr_k}{dx} + \frac{1}{4}\left(\theta_0\frac{dBr_i}{dx}Br_k + \frac{d\theta_0}{dx}Br_i\, Br_k\right) + \frac{1}{4}\sum_{j=2}^{n}a_j\left(Br_i\frac{dBr_j}{dx}Br_k\right)\right]dx\right]a_i$$

$$= \int_0^1\left[\left(16 + (2x+1)^3\right)Br_k + \frac{d\theta_0}{dx}\frac{dBr_k}{dx} - \frac{1}{4}\theta_0\frac{d\theta_0}{dx}Br_k\right]dx \qquad k = 2,3,\ldots,n \tag{22}$$

The above equations (22) are equivalent to the matrix form

$$(D + C)\,A = G \tag{23}$$

where the elements of the matrix $A, C, D$ and $G$ are $a_i, c_{i,k}, d_{i,k}$ and $g_k$, respectively, given by

$$d_{i,k} = \int_0^1\left[-\frac{dBr_i}{dx}\frac{dBr_k}{dx} + \frac{1}{4}\left(\theta_0\frac{dBr_i}{dx}Br_k + \frac{d\theta_0}{dx}Br_i\, Br_k\right)\right]dx \tag{24a}$$

$$c_{i,k} = \frac{1}{4}\sum_{j=2}^{n}a_j\int_0^1\left(Br_i\frac{dBr_j}{dx}Br_k\right)dx \tag{24b}$$

$$g_k = \int_0^1\left[\left(16 + (2x+1)^3\right)Br_k + \frac{d\theta_0}{dx}\frac{dBr_k}{dx} - \frac{1}{4}\theta_0\frac{d\theta_0}{dx}Br_k\right]dx \tag{24c}$$

The initial values of these coefficients $a_i$ are obtained by applying the Galerkin method to the BVP neglecting the nonlinear term in (19*a*). That is, to find initial coefficients, we will solve the system

$$DA = G \tag{25a}$$

where the matrices are constructed from

$$d_{i,k} = -\int_0^1 \frac{dBr_i}{dx}\frac{dBr_k}{dx}\,dx\,, \qquad \text{and} \qquad g_k = \int_0^1 \left[ \left(16 + (2x+1)^3\right) Br_k + \frac{d\theta_0}{dx}\frac{dBr_k}{dx} \right] dx \tag{25b}$$

Once the initial values of the parameters $a_i$ are obtained from Equation (25*a*), they are substituted into Equation (23) to obtain new estimates for the values of $a_i$. This iteration process continues until the converged values of the unknowns are obtained. Substituting the final values of the coefficients in Equation (20), we obtain an approximate solution of the BVP (19), and if we replace $x$ by $(x-1)/2$ in this solution we will obtain the approximate solution of the given BVP (18).

Using first 10 and 15 Bernoulli polynomials with 8 iterations, the absolute differences between exact and the approximate solutions are sown in **Table 5**. It is observed that the accuracy is found of the order nearly $10^{-6}$ and $10^{-8}$ on using 10 and 15 Bernoulli polynomials, respectively.

**Table 5. Approximate solutions of example 5 using 8 iterations.**

| $x$ | Approximate | Absolute Error | Approximate | Absolute Error |
|---|---|---|---|---|
| | Bernoulli polynomials = 10 | | Bernoulli polynomials = 15 | |
| 1.0 | 17.0000000000 | 0.0000000000 | 17.0000000000 | 0.0000000000 |
| 1.1 | 15.7554441298 | $1.041563\times10^{-5}$ | 15.7554545265 | $1.894190\times10^{-8}$ |
| 1.2 | 14.7733332492 | $8.415722\times10^{-8}$ | 14.7733333734 | $4.008605\times10^{-8}$ |
| 1.3 | 13.9977064922 | $1.418450\times10^{-5}$ | 13.9976922315 | $7.623683\times10^{-8}$ |
| 1.4 | 13.3885732769 | $1.848369\times10^{-6}$ | 13.3885714535 | $2.496207\times10^{-8}$ |
| 1.5 | 12.9166531985 | $1.346816\times10^{-5}$ | 12.9166667280 | $6.132234\times10^{-8}$ |
| 1.6 | 12.5599891326 | $1.086740\times10^{-5}$ | 12.5599999558 | $4.420697\times10^{-8}$ |
| 1.7 | 12.3017691184 | $4.412558\times10^{-6}$ | 12.3017646447 | $6.122859\times10^{-8}$ |
| 1.8 | 12.1289033378 | $1.444889\times10^{-5}$ | 12.1288889220 | $3.310400\times10^{-8}$ |
| 1.9 | 12.0310620046 | $9.373023\times10^{-6}$ | 12.0310526888 | $5.724404\times10^{-8}$ |
| 2.0 | 11.9999952268 | $4.773188\times10^{-6}$ | 11.9999999739 | $2.605950\times10^{-8}$ |
| 2.1 | 12.0290338956 | $1.372349\times10^{-5}$ | 12.0290475563 | $6.270952\times10^{-8}$ |
| 2.2 | 12.1127183984 | $8.874283\times10^{-6}$ | 12.1127272862 | $1.346861\times10^{-8}$ |
| 2.3 | 12.2465264382 | $4.699113\times10^{-6}$ | 12.2465218029 | $6.380614\times10^{-8}$ |
| 2.4 | 12.4266794677 | $1.280107\times10^{-5}$ | 12.4266666577 | $9.009668\times10^{-9}$ |
| 2.5 | 12.6500062226 | $6.222611\times10^{-6}$ | 12.6499999354 | $6.458229\times10^{-8}$ |
| 2.6 | 12.9138385465 | $7.607374\times10^{-6}$ | 12.9138461739 | $2.007236\times10^{-8}$ |
| 2.7 | 13.2159161538 | $9.772167\times10^{-6}$ | 13.2159259793 | $5.335987\times10^{-8}$ |
| 2.8 | 13.5542901664 | $4.452108\times10^{-6}$ | 13.5542856601 | $5.415501\times10^{-8}$ |
| 2.9 | 13.9272471869 | $5.807598\times10^{-6}$ | 13.9272414119 | $3.262486\times10^{-8}$ |
| 3.0 | 14.3333333333 | 0.0000000000 | 14.3333333333 | 0.0000000000 |

**Example 6.** **C**onsider a nonlinear differential equation [9] with the *Robin* boundary conditions [17]:

$$\frac{d^2u}{dx^2} = \frac{1}{2}(1+x+u)^3, \quad 0 < x < 1. \tag{26a}$$

$$u'(0) - u(0) = -1/2 \quad \text{and} \quad u'(1) + u(1) = 1. \tag{26b}$$

The exact solution of the problem is given by

$$u(x)=\frac{2}{2-x}-x-1\,.$$

In this case, solving the nonlinear BVP (26) by Modified Galerkin method, the approximate solution is assumed by

$$u(x)=\sum_{i=0}^{n}a_i Br_i(x),\qquad n\geq 1, \tag{27}$$

Now following the procedures described as in example-5 and with minor simplifications, the Equation (22) leads us

$$\sum_{i=0}^{n}\left[\int_0^1\left[-\frac{dBr_i}{dx}\frac{dBr_k}{dx}-\frac{3}{2}(1+x)^2 Br_i Br_k-\frac{3}{2}\sum_{j=0}^{n}a_j\,\{(1+x)Br_i Br_j Br_k\}\right.\right.$$
$$\left.\left.-\frac{1}{2}\sum_{j=0}^{n}a_j\left(\sum_{l-2}^{n}a_l\,(Br_i Br_j Br_l Br_k)\right)\right]dx-\{Br_i(1)Br_k(1)+Br_i(0)Br_k(0)\}\right]a_i$$
$$=\frac{1}{2}\int_0^1(1+x)^3 Br_k dx-\frac{1}{2}Br_k(0)-Br_k(1)\qquad\qquad k=0,1,2,\ldots,n \tag{28}$$

which can be written in a matrix form, similar to the system (23),

$$(D+C+B)\,A=G \tag{29a}$$

where the elements of $A,B,C,D$ and $G$ are $a_i$, $b_{i,k}$, $c_{i,k}$, $d_{i,k}$ and $g_k$, respectively, given by

$$d_{i,k}=\int_0^1\left(\frac{dBr_i}{dx}\frac{dBr_k}{dx}+\frac{3}{2}(1+x)^2 Br_i Br_K\right)dx+Br_i(1)Br_k(1)+Br_i(0)Br_k(0) \tag{29b}$$

$$c_{i,k}=\frac{3}{2}\sum_{j=0}^{n}a_j\int_0^1\{(1+x)Br_i Br_j Br_k\}\,dx \tag{29c}$$

$$b_{i,k}=\frac{1}{2}\sum_{j=0}^{n}a_j\left[\sum_{l-0}^{n}a_l\int_0^1(Br_i Br_j Br_l Br_k)dx\right] \tag{29d}$$

$$g_k=-\frac{1}{2}\int_0^1(1+x)^3 Br_k dx+\frac{1}{2}Br_k(0)+Br_k(1) \tag{29e}$$

To find the initial coefficients, on neglecting the nonlinear terms as we have done in example 5, we solve the reduced system

$$DA=G, \tag{30a}$$

where the elements of $D$, $A$ and $G$, respectively, are now

$$d_{i,k}=\int_0^1\left(\frac{dBr_i}{dx}\frac{dBr_k}{dx}+\frac{3}{2}(1+x)^2 Br_i Br_K\right)dx+Br_i(1)Br_k(1)+Br_i(0)Br_k(0) \tag{30b}$$

$$g_k=-\frac{1}{2}\int_0^1(1+x)^3 Br_k dx+\frac{1}{2}Br_k(0)+Br_k(1) \tag{30c}$$

**Table 6. Approximate solutions of examples using 8 iterations.**

| $x$ | Approximate | Absolute Error | Approximate | Absolute Error |
|---|---|---|---|---|
| | Bernoulli polynomials, 8 | | Bernoulli polynomials, 10 | |
| 0.0 | 00.0000000000 | 0.0000000000 | 00.0000000000 | 0.0000000000 |
| 0.1 | -0.0473680463 | $3.747086\times10^{-7}$ | -0.0473684172 | $3.803934\times10^{-9}$ |
| 0.2 | -0.0888891811 | $2.922364\times10^{-7}$ | -0.0888888959 | $6.992957\times10^{-9}$ |
| 0.3 | -0.1235296394 | $2.276838\times10^{-7}$ | -0.1235293980 | $1.372571\times10^{-8}$ |
| 0.4 | -0.1499995043 | $4.957444\times10^{-7}$ | -0.1500000116 | $1.157174\times10^{-8}$ |
| 0.5 | -0.1666665743 | $9.237249\times10^{-8}$ | -0.1666666693 | $2.681433\times10^{-9}$ |
| 0.6 | -0.1714291237 | $5.522990\times10^{-7}$ | -0.1714285563 | $1.508438\times10^{-8}$ |
| 0.7 | -0.1615383671 | $9.448797\times10^{-8}$ | -0.1615384751 | $1.360781\times10^{-8}$ |
| 0.8 | -0.1333328804 | $4.528913\times10^{-7}$ | -0.1333333287 | $4.621790\times10^{-9}$ |
| 0.9 | -0.0818186682 | $4.863450\times10^{-7}$ | -0.0818181840 | $2.134837\times10^{-9}$ |
| 1.0 | 00.0000000000 | 0.0000000000 | 00.0000000000 | 0.0000000000 |

The results are summarized in the **Table 6** that obtained on using 8 and 10 Bernoulli polynomials with 8 iterations at various points of the domain of the problems. It is observed that the approximate results converge monotonically to the exact solutions.

## 5. Conclusions

We have discussed, in details, the formulation of one dimensional linear and nonlinear second order boundary value problems by Galerkin weighted residual method, using Bernoulli polynomials which have been used as the trial functions in the approximation. Some numerical examples are tested. All the mathematical formulations and numerical computations have been evaluated by *MATHEMATICA* code. The computed solutions are compared with the exact solutions, and we have found a good agreement with the exact solution.